\documentclass{article} \usepackage{amssymb,amsmath,eucal}

\begin{document}

\def\R {{\mathbb R }}
 \def\C {{\mathbb C }}
  \def\Z{{\mathbb Z}} 
  \def\H{{\mathbb H}}

\def\vr{\mathbf {wr}}

\newcommand{\sgn}{\mathop{\mathrm{sgn}}\nolimits}

\def\SL{{\rm SL}}

\def\U{{\rm U}}
\def\O{{\rm O}} 
\def\Sp{{\rm Sp}} 
\def\SO{{\rm SO}}

\def\frm{\frak m}

\def\ov{\overline} 
\def\phi{\varphi} 
\def\epsilon{\varepsilon}
\def\kappa{\varkappa}

\def\le{\leqslant} 
\def\ge{\geqslant}

\def\frsl{\mathfrak{sl}_2}

\def\sl{\mathrm{SL}(2,\R)}
\def\slsim{\mathrm{SL}(2,\R)^\sim}

\renewcommand{\Re}{\mathop{\rm Re}\nolimits} 
\renewcommand{\Im}{\mathop{\rm Im}\nolimits}

\newcommand{\Arg}{\mathop{\rm Arg}\nolimits}








\def\cL{\mathcal L} 
\def\cM{\mathcal M}
\def\cH{\mathcal H}
\def\cK{\mathcal K}
\def\cE{\mathcal E}

\def\mn{\mathsf{M}}

\def\wt{\widetilde}

\begin{center}
{\Large\bf 
 Perelomov  problem and inversion of the Segal-Bargmann transform
}

\large\sc Neretin Yuri A.%
\footnote{Suppoted by the grant NWO-RFBR 047.011.2004.059}  

\end{center}

{\small We reconstruct a function by values of
its Segal-Bargmann transform at points of a lattice.}

{\bf 1. Formulation of the result.}
Fix $\tau>0$.
For a function $f\in L^2(\R)$, we define the coefficients
$$
\gamma_{m,k}=\int_{-\infty}^\infty e^{-ikx-\tau m x} f(x) e^{-x^2/4}
\,dx
$$
where $m$, $k$ range in $\Z$.
We intend to reconstruct $f$ by $\gamma_{m,k}$.
As Perelomov showed, this is impossible for $\tau>\pi$;
for $\tau\le\pi$, the problem is overdetermined
(see \cite{Per1}-\cite{Per2}, \cite{BiK}, more recent results
in \cite{Lyu}, \cite{Gro}).
There are
many ways for reconstruction of $f$. We propose a formula
(for $\tau<\pi$ ) that seems relatively simple and relatively closed. 

Denote $q:=e^{-2\pi\tau}$. Define the coefficicients
\begin{equation}
\cE_m(\tau)=
\frac{(-1)^{m} q^{m(m-1)/2}}{\prod_{l=1}^{\infty}(1-q^l)^3}
\sum_{j\ge 0} (-1)^j q^{j(j+2m+1)/2}
\label{cE}
\end{equation}

Then
$$
f(x)=\frac 1{2\pi}e^{x^2/4}
\sum_m \Bigl\{ \cE_m(\tau) e^{m\tau x} \sum_k \gamma_{m,k}e^{ikx}\Bigr\}
$$

The interior sum is an $L^2$-sum of a Fourier series,
the exterior sum is a.s. convergent series. 

\smallskip

{\sc Remark.} 1.  Our problem also is known in
the theory of recognition and separation
of waves (i.e., sound or electromagnetic oscillations).
A.J.E.M.Janssen \cite{Jan}  proposed several non-equivalent
formulae for reconstruction of
of the function $f$; the formula (\ref{cE})  easily follows
from his considerations. Hence, this note
clarifies Janssen's results and present them in a final
nice form.

2. Two-step way of reconstruction of $f$ was proposed by Yu.Lyubarskii,
he uses the Lagrange formula to interpolate the Segal--Bargman transform 
of a function by 
its values at points of the lattice. Then we apply the inversion
formula for the Segal--Bargmann transform.
I do not know, is it possible
 to produce a nice one-step
formula from this algorithm.
Our calculation is based on the same idea but
it gives another final result. 
 
\smallskip 

{\bf 2. Preliminaries on $\theta$-functions.}
Let $0<q<1$.
Denote
$$
\Theta(z;q):=(1-z)\prod_{n=1}^\infty (1-q^n)(1-z q^n)(1-z^{-1}q^n)=
\sum_{-\infty}^\infty (-1)^n z^n q^{n(n-1)/2}
$$
(this is the Jacobi triple identity, see, for instance,
\cite{Ahi}).
Obviously, 
$$\Theta(qz;q)=-z^{-1} \Theta(z;q)$$
Iterating this identity, we obtain
\begin{equation}
\Theta(q^n z;q)=(-z)^{-n} q^{-n(n-1)/2} \Theta(z;q)
\label{periodicity}
\end{equation}
The function 
\begin{equation}
\eta(z)=\exp\bigl\{-\frac 1{2\ln q}\ln^2|z|+\frac 12 \ln q\ln|z|\bigr\}
\label{asympt-1}
\end{equation}
satisfies the  requrence equation $\eta(qz)=|z|^{-1} \eta(z)$.
 Hence $|\Theta(z;q)|$ can be represented
in the form
\begin{equation}
|\Theta(z;q)|=\eta(z)\psi(z);\qquad \text{where $\psi(qz)=\psi(z)$}
\label{asympt-2}
\end{equation}

Obviously
$$
\Theta'(1;q)= \frac d{dx} \Theta(x;q)\Bigr|_{x=1} =-\prod(1-q^n)^3
$$
Differentiating (\ref{periodicity}) and subsituting $z=1$,
we obtain
\begin{equation}
\Theta'(q^n;q)=(-1)^{n+1} q^{-n(n-1)/2} \Theta'(1;q)
\label{shtrih}
\end{equation}

{\bf 3. Interpolation problem.}
Denote $g(x)=\frac1{2\pi}f(x)e^{-x^2/4}$.
Applying the Poisson summation formula
to the function $g(x)e^{-\tau mx}$, we obtain
$$
e^{m\tau x}\sum_{k=-\infty}^\infty \gamma_{m,k} e^{ikx}
=
\sum_{j=-\infty}^\infty g(x+2\pi j)e^{-2\pi\tau m j}
$$
Denote the right-hand side of this identity by $A_m$
Consider the function
$$
G_x(z):= \sum_{j=-\infty}^\infty g(x+2\pi j) z^j
$$ 
defined in the domain $\C\setminus 0$,
$$
G_x(q^m)=A_m
$$
We obtain
 an interpolation problem for holomorphic functions,
and solve it in a standard way (see \cite{Lev}).

\smallskip

Denote
\begin{equation}
\wt G_x(z)=
\sum_{n=-\infty}^\infty A_n \frac{\Theta(z;q)}
{(z-q^n)\Theta'(q^n;q)}
=
\sum_{n=-\infty}^\infty 
A_n\frac{(-1)^{n+1} q^{n(n-1)/2}}
{\prod(1-q^j)^3} 
\frac{\Theta(z;q)}{(z-q^n)}
\label{wt-G}
\end{equation}

Obviously,
\begin{equation}
G_x(q^n)=\wt G_x(q^n)
\end{equation}
Hence,
\begin{equation}
 G_x(z)=\wt G_x(z) + \Theta(z;q)\alpha(z)
\label{identity}
\end{equation}
for certain function $\alpha(z)$
holomorphic in $\C\setminus 0$.

\smallskip

{\sc Lemma.}  $G_x(z)=\wt G_x(z)$, i.e., $\alpha(z)=0$.

\smallskip

Our final formula is a corollary of this lemma.
Indeed, $g(x)$ is the Laurent coefficient
of $G_x(z)$ in $z^0$; it remains to evaluate
the Laurent expansion of 
$$(z-q^n)^{-1}\Theta(z;q)= 
(z-q^n)^{-1}\sum_{l=-\infty}^\infty (-1)^l z^l q^{l(l-1)/2}$$
Assuming $|z|>q^n$, we obtain
$$(z^{-1}+z^{-2}q^n+z^{-3} q^{2n}+\dots)\cdot \sum_{l=-\infty}^\infty (-1)^l z^l q^{l(l-1)/2}$$
and we obtain (\ref{cE}) as a coefficient in the front of
$z^0$.

{\bf 4. Proof of Lemma.}
We represent the identity (\ref{identity}) in the form
\begin{equation}
 G_x(z)/ \Theta(z;q) =\wt G_x(z)/ \Theta(z;q)  + \alpha(z)
\label{identity-2}
\end{equation}
For a function $\Phi(z)$ we denote 
$$\cM_k[\Phi]
:=\max_{|z|=q^{k+1/2}} 
|\Phi(z)|
$$
We intend to analize the behavior
of these maxima for summands of (\ref{identity-2})
  as $k\to\pm\infty$.

A) First,
$$
\infty>\int_\R |f(x)|^2 dx
=\int_0^{2\pi}\Bigl( \sum_{j=-\infty}^\infty |f(x+2\pi j)|^2
\Bigr)\,dx
$$
Hence (by the Fubini theorem) the value
$$V_x:=\sum_{j=-\infty}^\infty |f(x+2\pi j)|^2$$ 
is finite for almost
all $x$. 

\smallskip

B) By the Schwartz inequality,
\begin{multline*}
|G_x(z)|=\Bigl|\sum f(x+2\pi j) e^{-(x+2\pi j)^2/4} z^j\Bigr|
\le\\ \le
\Bigl( \sum |f(x+2\pi j)|^2\Bigr)^{1/2}
\Bigl(\sum e^{-(x+2\pi j)^2/2} |z|^{2j}\Bigr)^{1/2}
=\\
=V_x^{1/2}\cdot\Bigl[ e^{-x^2} \Theta( -|z|^2 e^{-2\pi x-2\pi^2} ;
e^{-4\pi^2})\Bigr]^{1/2}
\end{multline*}
Applying (\ref{asympt-1})-(\ref{asympt-2}),
 we obtain for $|G_x(z)|$
 an upper estimate of the form
\begin{equation}
|G_x(z)|\le \exp\bigl\{\frac 1{4\pi^2}\ln^2|z|
+ 
O(\ln|z|)+ O(1)\bigr\}
\label{gx}
\end{equation}
In particular,
$$
|A_m|=|G_x(q^m)|\le 
\exp\bigl\{ \frac {\ln^2 q}{4\pi^2} \,\,  m^2+ O(m)+O(1)
\bigr\}
$$

By (\ref{shtrih}),
$$\Theta'(q^m;q)=\exp\bigl\{-m^2 \ln q/2+O(m)+O(1)\bigr\}$$
Since $(-\ln q)=2\pi\tau<2\pi^2$,
 we obtain the following estimate
$$
|A_m / \Theta'(q^m;q)|\le \exp\bigl\{- \epsilon m^2\bigr\}
$$

C)
Consider the summand $\wt G_x(z)/\Theta(z;q)$  in  (\ref{identity-2})
\begin{multline*}
\cM_k\bigl[\wt G_x(z)/\Theta(z;q)]=
\cM_k\Bigl[
\sum_m\frac{A_m }{ \Theta'(q^m;q)}\cdot \frac 1{z-q^m}\Bigr]
\le \sum \frac{e^{-\epsilon m^2}} {|q^{k+1/2}-q^m|}
\end{multline*}
Next,
$$
|q^{k+1/2}-q^m|=q^m |1-q^{-m+k+1/2}|\ge q^m (1-q^{1/2})
$$
This implies the boundedness of the sequence $\cM_k[\cdot]$.

Secondly,
$$
|q^{k+1/2}-q^m|\ge q^{k+1} (1-q^{1/2})
$$
Hence,
$\cM_k[\cdot]$ tends to 0 as $k\to-\infty$.

D) By (\ref{asympt-2})
$$
\cM_k\bigl[
\Theta(z)^{-1}
\bigr]\sim \eta(z)^{-1}\Bigr|_{|z|=q^{k+1/2}}
$$
By (\ref{asympt-1}), (\ref{gx})
$$
\cM_k\bigl[G_x(z)/\Theta(z;q)\bigr]\to 0\qquad \text{as $k\to\pm\infty$}
$$

E) 
We have 
$$\cM_k\bigl[ \alpha(z)\bigr]\le 
\cM_k\bigl[ G_x(z) / \Theta(z;q)\bigr] 
+ \cM_k\bigl[ \wt G_x(z) / \Theta(z;q)\bigr] 
$$ 
Thus
$\cM_k[ \alpha(z)]$
 tends to 0 as $k\to-\infty$;
 and remains bounded as  $k\to +\infty$.
Since $\alpha(z)$ is holomorphic in $\C\setminus 0$,
we have $\alpha(z)=0$.

{\bf Acknowledgements.}
I am grateful to Yu. Lyubarskii who explained 
me various tricks related to entire functions.
I  also thank H.Feichtinger, A.M.Vershik,
and K.Grochenig for discussion of the subject.

{\sf Math.Phys. Group,
Institute of Theoretical and Experimental Physics,

B.Cheremushkinskaya, 25, Moscow 117259

\& University of Vienna, Math. Dept.,
Nordbergstrasse, 15, Vienna 1090, Austria}

neretin@mccme.ru


\begin{thebibliography}{cc}
\bibitem{Ahi}
Akhiezer, N. I.
{\it Elements of the theory of elliptic functions.}
Translations of Mathematical Monographs, 79.
American Mathematical Society, Providence, RI, 1990. 



\bibitem{BiK}
Bargmann, V.; Butera, P.; Girardello, L.; Klauder, J. R.
{\it On the completeness of the coherent states.}
Rep. Mathematical Phys. 2 1971 no. 4, 221--228.

\bibitem{Gro}
Gr\"ochenig, K.
{\it Foundations of time-frequency analysis.}
Birkh\"auser, 2001


\bibitem{Jan}
Janssen A.J.E.M., {\it Signal analytic proof
of two basic results on lattice expansion,}
 Applied and computational harmonic analysis,
I (1994), 330-354

\bibitem{Lev}
Levin B.Ya. {\it Distribution of zeros of entire functions.}
 American Mathematical Society, Providence, R.I. 1964 
 
 
 \bibitem{Lyu}
 Lyubarskii, Yu. I. {\it 
  Frames in the Bargmann space of entire functions.}
    Entire and subharmonic functions,  167--180, Adv. Soviet Math., 11, 
    Amer. Math. Soc., Providence, RI, 1992.
 
\bibitem{Per1}
Perelomov, A. M. {\it Remark on the completeness of the coherent state system.}   Teoret. Mat. Fiz.  6  (1971), no. 2, 213--224.

\bibitem{Per2}
Perelomov, A. {\it Generalized
 coherent states and their applications.} 
  Springer-Verlag, Berlin, 1986. xii+320 pp


 \end{thebibliography}
\end{document}